\documentclass[12pt,reqno]{amsart}
\usepackage{amssymb,delarray}
\usepackage{amsfonts}
\usepackage{epsfig}
\usepackage[all]{xy}
\usepackage{amscd}
\usepackage{epigraph}

\makeindex{}

\newtheorem{thm}{Theorem}
\newtheorem{lem}{Lemma}
\newtheorem{prop}{Proposition}
\newtheorem{claim}{Claim}

\newtheorem{cor}{Corollary}

{\catcode`\@=11
\gdef\n@te#1#2{\leavevmode\vadjust{%
 {\setbox\z@\hbox to\z@{\strut#1}%
  \setbox\z@\hbox{\raise\dp\strutbox\box\z@}\ht\z@=\z@\dp\z@=\z@%
  #2\box\z@}}}
\gdef\leftnote#1{\n@te{\hss#1\quad}{}}
\gdef\rightnote#1{\n@te{\quad\kern-\leftskip#1\hss}{\moveright\hsize}}
\gdef\?{\FN@\qumark}
\gdef\qumark{\ifx\next"\DN@"##1"{\leftnote{\rm##1}}\else
 \DN@{\leftnote{\rm??}}\fi{\rm??}\next@}}

\begin{document}

\baselineskip=14.pt plus 2pt 

\title[Galois groups of dualizing coverings]
{On the Galois groups of the dualizing coverings for plane curves}
\author[Vik.S.~Kulikov]{Vik.S. Kulikov}

\address{Steklov Mathematical Institute}
\email{kulikov@mi.ras.ru}

\thanks{This
research was partially supported by grants of NSh-2998.2014.1, RFBR
14-01-00160, and by AG Laboratory HSE, RF government grant, ag.
11.G34.31.0023.. }


\begin{abstract}
Let $C_1$ be an irreducible component of a reduced projective curve
$C\subset \mathbb P^2$ defined over the field $\mathbb C$, $\deg
C_1\geq 2$, and let $T$ be the set of lines $l\subset \mathbb P^2$
meeting $C$ transversally. In the article, we prove that for a line
$l_0\in T$ and any two points $P_1,P_2\in C_1\cap l_0$ there is a
loop $l_t\subset T$, $t\in [0,1]$, such that the movement of the
line $l_0$ along the loop $l_t$ induces the transposition of the
points $P_1$, $P_2$ and the identity permutation of the other points
of $C\cap l_0$.
\end{abstract}

\maketitle

\setcounter{tocdepth}{2}


\def\st{{\sf st}}

\section*{Introduction}

Let $C_i\subset \mathbb P^2$, $1\leq i\leq k$, be  irreducible
reduced curves defined over the field $\mathbb C$, $\deg C_i=d_i\geq
2$, and $C=C_1\cup\dots\cup C_k$, $d=\deg C=d_1+\dots +d_k$. Denote
by $\nu:\overline C\to C$ the normalization of the curve $C$ and
consider a point $p\in \overline C$ and its image $P=\nu(p)\in
\mathbb P^2$. Choose homogeneous coordinates $(x_1,x_2,x_3)$ in
$\mathbb P^2$ such that $P=(0,0,1)$. We can choose a local parameter
$t$ in a complex analytic neighborhood $U\subset \overline C$ of the
point $p$ such that  the regular map $\nu$ is given by
\begin{equation} \label{nu} x_1=\sum_{i=s_p}^{\infty} a_it^i, \, \,
x_2=t^{r_p}, \, \, x_3=1, \end{equation} where $a_{s_p}\neq 0$  and
$s_p>r_p\geq 1$. The integer $r_p$ is called the {\it multiplicity}
of the germ $\nu(U)$ of the curve $C$ at $P=\nu(p)$, the line
$l_p=\{ x_1=0\}$ is called a {\it tangent line} to $C$ at $P$, and
the integer $s_p$ is called the {\it tangent multiplicity} of the
germ $\nu(U)$ at $P$.

Let $\hat C\subset \hat{\mathbb P}^2$ be the dual curve to the curve
$C$ (the curve $\hat C$ consists of the tangents $l_p$, $p\in
\overline C$, to $C$). The graph of the correspondence between $C$
and $\hat C$ is a curve $\check C$ (the so called Nash blow-up of
$C$) in $\mathbb P^2\times \hat{\mathbb P}^2$ which lies in the
incidence variety $I=\{ (P,l)\in \mathbb P^2\times \hat{\mathbb
P}^2\mid P\in l\}$,
$$\check C=\{ (\nu(p),l_p)\in I\mid p\in \overline C\, \, \text{and}\, \, l_p \, \,
\text{is the tangent line to}\, \, C\,\, \text{at}\,\, \nu(p)\in C\}.$$

In the sequel, by $L_p\subset \hat{\mathbb P}^2$, we denote the line
dual to the point $\nu(p)\in C\subset \mathbb P^2$.

Let $\text{pr}_1:\mathbb P^2\times \hat{\mathbb P}^2\to \mathbb P^2$
and $\text{pr}_2:\mathbb P^2\times \hat{\mathbb P}^2\to \hat{\mathbb
P}^2$ be the projections to the factors, $X=\text{pr}_1^{-1}(C)\cap
I$, and $f':X\to \hat{\mathbb P}^2$ the restriction of $\text{pr}_2$
to $X$. Obviously, $f'^{-1}(l)$ consists of the points
$(P,l)\in\mathbb P^2\times \hat{\mathbb P}^2$ such that $P\in C\cap
l$ and hence $\deg f'=\deg C=d$.

Denote by $\nu':Z\to X$ the normalization of $X$ and by $f=f'\circ
\nu':Z\to \hat{\mathbb P}^2$. We have $\deg f=d$. We call $f$ the
{\it dualizing covering} for $C\subset \mathbb P^2$. Obviously, the
variety $Z$ is isomorphic to the fibre product $\overline C\times_C
X$ of the normalization $\nu:\overline C\to C$ and the projection
$\text{pr}_1:X\to C$. The projection $\text{pr}_1:X\to C$ gives on
$X$ a structure of a ruled surface and it induces a ruled structure
on $Z$ over the curve $\overline C$, $\rho :Z\to \overline C$,
$\rho^{-1}(p):=F_p\simeq \mathbb P^1$ for $p\in \overline C$. Note
that the curve $\widetilde C$ is a section of this ruled structure,
where $\widetilde {C}=\nu'^{-1}(\check C)\subset Z$, and the image
$f(F_p)$ of a fibre $F_p$ is the line $L_p\subset \hat{\mathbb P}^2$
dual to the point $\nu(p)\in C\subset \mathbb P^2$.

Denote by $\overline B\subset \hat{\mathbb P}^2$ the branch locus of
$f$, choose a point $l_0\in \hat{\mathbb P}^2\setminus \overline B$,
and number the points of $f^{-1}(l_0)$. In this case the covering
$f$ induces a homomorphism $f_*:\pi_1(\hat{\mathbb P}^2\setminus
\overline B,l_0)\to \Sigma_d$ from the fundamental group
$\pi_1(\hat{\mathbb P}^2\setminus \overline B,l_0)$ to the symmetric
group $\Sigma_d$ acting on the fibre $f^{-1}(l_0)$. The image, $Im
f_*:=G\subset \Sigma_d$, is called the {\it Galois group} of the
covering $f$.

\begin{thm} \label{main0} Let $f:Z\to \hat{\mathbb P}^2$ be  the
dualizing covering for a reduced curve $C\subset \mathbb P^2$, $\deg
C=d$, and $C_1,\dots, C_k$ the irreducible components of $C$, $\deg
C_i=d_i\geq 2$. Then the Galois group of $f$ is  $G\simeq
\Sigma_{d_1}\times \dots\times \Sigma_{d_k}$.
\end{thm}

The following theorem describes properties of dualizing coverings.

\begin{thm} \label{main1} Let $C$ be as in Theorem {\rm \ref{main0}}
and $f:Z\to \hat{\mathbb P}^2$ the dualizing covering for $C$. Then
$Z$ is a non-singular surface consisting of $k$ irreducible
components and $f$ is a degree $d$ finite covering.

The branch locus of $f$ is $\overline B=\hat C\cup \hat L$, where
$\hat L=\displaystyle \bigcup_{r_p\geq 2} L_p$, $L_p$ are the lines
dual to the points $\nu(p)\in C$ and the union is taken over all
$p\in \overline C$ for which the multiplicity $r_p\geq 2$.

The ramification locus of $f$ is $\overline R=\widetilde  C\cup
\widetilde F$, where $\widetilde F=\displaystyle \bigcup_{r_p\geq 2}
F_p$ and the union is taken over all $p\in \overline C$ for which
$r_p\geq 2$.

The local degree $\deg_q f$ of $f$ at a point $q=F_p\cap \widetilde
C$ is equal to the tangent multiplicity $s_p$, and $\deg_q f=r_p$ at
all points $q\in F_p\setminus \widetilde C$. For all points
$q\in\widetilde C\setminus \widetilde F$ the local degree $\deg_q
f=2$.
\end{thm}

For given reduced projective curve $C\subset \mathbb P^2$, $\deg
C=n$, let $T_C$ be the set of lines $l\subset \mathbb P^2$ meeting
$C$ transversally. Let $l_t\subset T_C$, $t\in [0,1]$, be a loop and
let $l_0\cap C=\{ P_1,\dots, P_n\}$. Then the movement of the line
$l_0$ along the loop $l_t$ defines $n$ paths $\psi_i(t)=l_t\cap
C\subset \mathbb P^2$, $i=1,\dots, n$, starting and ending at the
points $P_1,\dots, P_n$ and, consequently, induces a  permutation of
the points $P_1,\dots, P_n$ called the {\it monodromy} of the points
$P_1,\dots, P_n$ along the loop $l_t$ (the start point
$P_i=\psi_i(0)$ of the path $\psi_i(t)$ maps to the end point
$\psi_i(1)\in l_0\cap C$).

\begin{cor} \label{main2}
Let $C_1$, $\deg C_1\geq 2$, be an irreducible component of a
reduced curve $C\subset \mathbb P^2$. For a line $l_0\subset T_C$
and any two points $P_1,P_2\in C_1\cap l_0$ there is a loop
$l_t\subset T_C$, $t\in [0,1]$, such that the monodromy along the
loop $l_t$ is the transposition of the points $P_1$, $P_2$ and the
identity permutation of the other points in $C\cap l_0$.
\end{cor}

The proof of Theorems 1, 2 and Corollary 1 will be given in Section
\ref{thm1}. In Section \ref{sec1}, we give some remarks on the
actions of finite groups on finite sets and prove Proposition
\ref{G+J} which plays the crucial role in the computation of the
Galois groups of the dualizing coverings in the proof of Theorem
\ref{main0}. In Section \ref{cover}, we remind some properties of
finite ramified coverings and begin to prove Theorem \ref{main0}.

For background results on plane algebraic curves and  dual to them
we refer to \cite{Br} and \cite{Wall94}.

\section{Some remarks on the actions of finite groups}\label{sec1}
Let $I$ be a finite set and $\Sigma_I$ the symmetric group acting on
$I$. An embedding $I_1\subset I_2$ defines the natural embedding
$\Sigma_{I_1}\subset \Sigma_{I_2}$. In we sequel, we will assume
that each finite set $I$ is a subset of an integer segment $[1,d]=\{
1,2,\dots, d\}$, so that $\Sigma_I\subset \Sigma_{[1,d]}:=\Sigma_d$.

Let $J=\{I_{1},\dots, I_{k}\}$ be a partition of the segment $[1,d]$.
The partition $J$ defines an embedding of the group
$\Sigma_J:=\Sigma_{I_1}\times \dots\times \Sigma_{I_k}$ into $\Sigma_d$.

We say that a partition $J$ of $[1,d]$ is {\it invariant} under the
action of a subgroup $G\subset \Sigma_d$ if $g(I_j)=I_j$ for all
$g\in G$ and all $I_j\in J$.

Let $J_i=\{I_{i,1},\dots, I_{i,k_i}\}$, $i=1,2$, be two partitions
of $[1,d]$. We say that the partition $J_1$ is
{\it thinner} than $J_2$ (resp., $J_2$ is {\it thicker} than $J_1$)
and write $J_1\preccurlyeq J_2$ if for each $j$, $1\leq j\leq k_1$,
there is $t(j)$ such that $I_{1,j}\subset I_{2,t(j)}$. For any two
partitions $J_i=\{I_{i,1},\dots, I_{i,k_i}\}$, $i=1,2$, denote by
$J_1\oplus J_2$ the thinnest partition of $[1,d]$ such that
$J_1\preccurlyeq J_1\oplus J_2$ and $J_2\preccurlyeq J_1\oplus J_2$.

\begin{claim} \label{cl} Let $G$ be a subgroup of $\Sigma_d$ and let
$J_i=\{I_{i,1},\dots, I_{i,k_i}\}$, $i=1,2$, be two partitions of
$[1,d]$ such that $\Sigma_{J_i}\subset G$ for $i=1,2$. Then
$\Sigma_{J_1\oplus J_2}\subset G$.
\end{claim}
\proof Obvious. \qed \\

It follows from Claim \ref{cl} that for each subgroup $G$ of
$\Sigma_d$ there is the thickest partition of $[1,d]$ ({\it denote
it by} $J_G$) such that $\Sigma_{J_G}\subset G$.

Let $\sigma= c_1\cdot \, .\, .\, .\cdot c_n\in \Sigma_d$ be the
factorization of $\sigma$ into the product of cycles with disjoint
orbits. The number $n_{\sigma}=n$ will be called the {\it length} of
cycle factorization.

\begin{lem} \label{H+J} Let $H$ be a subgroup of $\Sigma_d$
generated by a set of transpositions and a permutation $\sigma$, and
let $\sigma= c_1\cdot \, .\, .\, .\cdot c_{n_{\sigma}}$ be the
factorization of $\sigma$ into the product of cycles with disjoint
orbits. Assume that for each $i$, $1\leq i\leq n_{\sigma}$, there is
a partition $J_i=\{ I_{i,1},\dots,I_{i,k_i}\}$ of $[1,d]$ invariant
under the action of $\sigma$ and such that
\begin{itemize}
\item[($i$)]
for each $I_{i,j}\in J_i$ there is at most one cycle $c_{m(i,j)}$
entering into the factorization of $\sigma$ such that the cycle
$c_{m(i,j)}$ acts non-trivially on $I_{i,j}$,
\item[($ii$)] the cycle $c_{i}$ acts non-trivially on $I_{i,1}$
and the length of the cycle $c_{i}$ is strictly less than the
cardinality of $I_{i,1}$,
\item[($iii$)] the group $H$ acts transitively on $I_{i,1}$.
\end{itemize} Then $H=\Sigma_{J_H}$ and, in particular, $H$
is generated by transpositions.
\end{lem}
\proof Consider a set $I_{i,1}$. Let $l_{i,1}$ be its cardinality
and let $l_i$ be the length of the cycle $c_i$. We have
$l_i<l_{i,1}$ and it follows from ($i$) and ($iii$) that there
exists a transposition $\tau \in H\cap \Sigma_{I_{i,1}}$ such that
it commutes with $c_j$ if $j\neq i$ and it transposes an element
entering in the cycle $c_i$ and an element of $I_{i,1}$ which does
not enter in $c_i$. Without loss of generality, we can assume that
$I_{i,1}=\{ 1, 2, \dots, l_{i,1}\}$, $c_{i}=(1,2,\dots, l_{i})$, and
$\tau=(l_i, l_i+1)$. Therefore
$$\sigma^{-j}\tau\sigma^{j}=(l_i-j,l_i+1)\in H$$
for $j=1,\dots,l_i-1$ and hence
$H\cap\Sigma_{l_i+1}=\Sigma_{l_i+1}$, since the subgroup of $H$,
generated by the transpositions $\sigma^{-j}\tau\sigma^{j}$,
$j=0,1,\dots,l_i-1$, acts transitively on the set $\{ 1,2,\dots,
l_i+1\}$. If we apply conditions ($i$) -- ($iii$) $l$ times, where
$l=l_{i,1}-l_i$, we obtain that $\Sigma_{I_{i,1}}\subset H$ for each
$i$ and hence $\sigma$ is a product of some transpositions belonging
to $H$. \qed

The following proposition is an easy consequence of Claim \ref{cl}
and Lemma \ref{H+J}.
\begin{prop} \label{G+J} Let $G$ be a subgroup of $\Sigma_d$
generated by some set of transpositions and by permutations
$\sigma_1,\dots , \sigma_m$. Assume that for each $i$, $i=1,\dots,
m$, there are partitions $J_{i,j}$ of $[1,d]$, $1\leq j\leq
n_{\sigma_i}$, such that the subgroup $H_i$ of $G$, generated by
transpositions and by $\sigma_i$, and the partitions  $J_{i,j}$
satisfy the conditions of Lemma {\rm \ref{H+J}}. Then
$G=\Sigma_{J_G}$.
\end{prop}

\begin{cor} \label{coro1} Let $G\subset \Sigma_d$ satisfy the conditions
of Proposition {\rm \ref{G+J}.} Assume that there is a partition
$J=\{ I_1,\dots, I_k\}$ of $[1,d]$ such that $G$ leaves invariant
the partition $J$ and acts transitively on each $I_j\in J$, $1\leq
j\leq k$. Then $J=J_G$ and $G=\Sigma_{J}$.
\end{cor}

\section{Coverings}\label{cover}
By a {\it covering}\, we understand a branched covering, that is a
finite morphism $f:Z\to Y$ from a normal projective surface $Z$ onto
a non-singular irreducible projective surface $Y$. To each covering
$f$ we associate the branch locus $\overline B\subset Y$, the
ramification locus $\overline R\subset f^{-1}(\overline B)\subset
Z$, and the unramified part $Z\setminus f^{-1}(\overline B)\to
Y\setminus \overline B$ (which is the maximal unramified
subcovering). As is usual for unramified coverings of degree $d$,
there is a homomorphism $f_*$ which acts from the fundamental group
$\pi_1( Y\setminus \overline B,p_0)$ to the symmetric group
$\Sigma_d$ acting on the points of $f^{-1}(p_0)$. The homomorphism
$f_*$ (called {\it monodromy} of $f$) is defined by $f$ uniquely if
we number the points of $f^{-1}(p_0)$; reciprocally, according to
Grauert-Remmert-Riemann-Stein Extension Theorem (see, for example,
\cite{St}) the conjugacy class of $f_*$ defines $f$ up to an
isomorphism. The image $G\subset \Sigma_d$ of $f_*$ is a transitive
subgroup of $\Sigma_d$ if $Z$ is irreducible and in general case the
number of connected components of $Z$ is equal to the number of
orbits of the action of $G$ on $f^{-1}(p_0)$.

An element $\gamma_q$, $q\in \overline B\setminus \text{Sing}\,
\overline B$, of the fundamental group $\pi _1(Y\setminus \overline
B, p_0)$ is called a {\it geometric generator} if it is represented
by a loop $\Gamma_q$ of the following form. To define $\Gamma_q$,
let $L \subset Y$ be a curve meeting $\overline B$ transversely at
$q$ and let $S^1 \subset L$ be a circle of small radius with center
at $q$. The choice of an orientation on $Y$ defines an orientation
on $S^1$. Then $\Gamma_q$ is a loop consisting of a path $l$ in
$Y\setminus \overline B$ joining $p_0$ with a point $q_1\in S^1$,
the loop $S^1$ (with positive direction) starting and ending at
$q_1$, and a return path to $p_0$ along $l$ in the opposite
direction (of course, we must note that a geometric generator
$\gamma_q$ depends not only on $q$, but it depends also on the
choice of the path $l$). Note that if $Y$ is simply connected then
$\pi _1(Y\setminus \overline B, p_0)$ is generated by geometric
generators.

In the sequel, we will assume that the covering $f$ satisfies some
additional conditions. The first of them is:
\begin{itemize}
\item[($R_0$)]
{\it If for an irreducible component $\overline B_i$ of the branch
curve $\overline B$ the image $f_*(\gamma_{q_i})$ of a geometric
generator $\gamma_{q_i}\in\pi_1(Y\setminus \overline B, p_0)$,
$q_i\in \overline B_i$, is not a transposition, then $\overline B_i$
is a smooth curve and $f_{|\overline R_{i,j}}:\overline R_{i,j}\to
\overline B_{i}$ is an isomorphism for all $j$, $1\leq j\leq n$,
where $\overline R\cap f^{-1}(\overline B_i)=\overline R_{i,1}\cup
\dots, \cup\overline R_{i,n}$ is the decomposition of $\overline
R\cap f^{-1}(\overline B_i)$ into the union of irreducible
components}.
\end{itemize}
Let $r_{i,j}$ be the ramification multiplicity of $f$ along
$\overline R_{i,j}$ (that is, the local degree of $f$ at a generic
point of $\overline R_{i,j}$), then the cycle type of the
permutation $f_*( \gamma_{q_i})\in \Sigma_d$ is $(r_{i,1},
\dots,r_{i,n})$.

Let for $\overline B_1$ the image $f_*(\gamma_{q_1})$ be not a
transposition and let $\overline R_{1},\dots, \overline R_{n}$ be
the irreducible components of $\overline R\cap f^{-1}(\overline
B_1)$. For each point $o\in \overline B_i\cap \text{Sing}\,
\overline B$ let us choose a very small (in complex analytic
topology) neighbourhood $W_o\subset Y$ of the point $o$. Denote by
$B_{1}:=\overline B_1\setminus (\displaystyle
\bigcup_{o\in\text{Sing}\, \overline B}\overline W_o)$ and
$R_{j}:=\overline R_{j}\cap f^{-1}(B_{1})$, where $\overline W_o$ is
the closure of $W_o$ in $Y$. The following Lemma is well known.

\begin{lem}\label{rami1}
There are neighbourhoods $U_1\subset Y$ and
$V_{j}\subset Z$, $j=1,\dots,n$, such that
\begin{itemize}
\item[($i$)] $U_1\cap\overline B=B_{1}$ and $V_{j}\cap \overline R_{j}=R_{j}$,
\item[($ii$)] $U_1$ is biholomorphic to $B_{1}\times D_1$
and $V_{j}$ is biholomorphic to $R_{j}\times D_2$,
where $D_1=\{ u_1\in \mathbb C\mid |u_1|<1\}$ and
$D_2=\{ v_1\in \mathbb C\mid |v_1|<1\}$
are discs in $\mathbb C$,
\item[($iii$)] the restriction of $f$ to each $V_{j}$
is proper and $f(V_{j})=U_1$,
\item[($vi$)] if $u_2$ is a local parameter on $B_{1}$
at a point $p\in B_{1}$ and
$v_{2,j}=f_{|R_{j}}^*(u_2)$ is the local parameter at the point
$q_{j}=f_{|R_{j}}^{-1}(p)$ on $R_{j}$, then $f:V_{j}\to U_1$ is
given by $u_1=v_1^{r_1}$ and $u_2=v_{2,j}$ in a neighbourhood of the
point $q_{j}$.
\end{itemize}
\end{lem}

Consider a neighbourhood $U_1$ the existence of which is claimed in
Lemma \ref{rami1}. Let $\text{pr}:U_1\to B_1$ be the projection
defined by bi-holomorphic isomorphism $U_1\simeq B_{1}\times D_1$.
Let us choose a point $q_1\in B_1$, a point $p_1\in
\text{pr}^{-1}(q_1)\simeq D_1$ lying in the circle $\delta=\{
q_1\}\times \{ u_1\in \mathbb C \mid |u_1|= \frac{1}{2}\}\subset \{
q_1\}\times D_1\subset U_1\setminus B_1$ (let, for definiteness,
$u_1=\frac{1}{2}$ at the point $p_1$), and a path $l_1 \subset
Y\setminus \overline B$ connecting the points $p_0$ and $p_1$. The
choice of $l_1$ defines homomorphisms
$\text{im}_{l_1}:\pi_1(U_1\setminus B_1,p_1)\to \pi_1(Y\setminus
\overline B,p_0)$ and $\varphi_{1}=f_*\circ
\text{im}_{l_1}:\pi_1(U_i\setminus B_1,p_1)\to \Sigma_d$. Denote by
$H_{B_1}$ the image $\varphi_{1}(\pi_1(U_1\setminus B_1,p_{1}))$ in
$G$. Let $\gamma_{q_1}\in \pi_1(U_1\setminus B_1,p_1)$ be a
geometric generator represented by the circle $\delta$ (the circuit
in positive direction). The cycle type of the permutation
$\sigma_1=\varphi_1( \gamma_{q_1})$ is $(r_{1,1}, \dots,r_{1,n})$.

Let a set $S=\{ o_1,\dots,o_m\}$ be the intersection of $\overline
B_1$ and $\text{Sing}\, \overline B$. For a point $o_i\in S$ we
choose a small (in complex analytic topology) simply connected
neighbourhood $U_{o_i}\subset Y$ of the point $o_i$ such that the
number $k_i$ of the connected components $V_{o_i,1},\dots
,V_{o_i,k_i}$ of $f^{-1}(U_{o_i})$ is equal to the number of points
belonging to $f^{-1}(o_i)$. In addition, $U_{o_i}$ can be chosen so
that $\overline W_{o_i}\subset U_{o_i}$, where $W_{o_i}$ is the
neighbourhood of $o_i$ used in the definition of the neighbourhood
$U_1$.

Let us choose points $q_{o_i}\in B_1\cap U_{o_i}$ and paths
$l'_{o_i}\subset B_1$ connecting, respectively, the point $q_{1}$
with the points $q_{o_i}$. Let $l_{o_i}=\{ p\in U_1\setminus B_1\mid
\text{pr}(p)\in l'_{o_i}, u_1(p)=\frac{1}{2}\}$ be paths connecting
the point $p_1$, respectively, with points $p_{o_i}=
\text{pr}^{-1}(q_{o_i})\cap l_{o_i}\in U_1\cap U_{o_i}$ (without
loss of generality, we can assume that $\text{pr}^{-1}(q)\subset
U_{o_i}$ if $q\in B_1\cap U_{o_i}$). Denote by $\widetilde l_{o_i}$
the composition of paths $l_1$ and $l_{o_i}$ connecting the point
$p_0$ with the point $p_{o_i}$.

The path $\widetilde l_{o_i}$ defines
homomorphisms $\text{im}_{\widetilde l_{o_i}}:\pi_1(U_{o_i}
\setminus (U_{o_i}\cap\overline B),p_{o_i})\to
\pi_1(Y\setminus \overline B,p_0)$ and
$\varphi_{\widetilde 1_{o_i}}=f_*\circ \text{im}_{\widetilde l_{o_i}}:
\pi_1(U_{o_i}\setminus (U_{o_i}\cap\overline B),p_{o_i})\to
\Sigma_d$. Denote by $H_{o_i}$ the image
$\varphi_{\widetilde 1_{o_i}}(
\pi_1(U_{o_i}\setminus (U_{o_i}\cap\overline B),p_{o_i}))$ in $G$.

Similarly, if $U=U_1\cup (\displaystyle \bigcup_{o_i\in S}U_{o_i})$,
then the path $l_1$ defines homomorphisms
$\text{im}_{l_1}:\pi_1(U\setminus (U\cap \overline B),p_1)\to
\pi_1(Y\setminus \overline B,p_0)$ and $\varphi=f_*\circ
\text{im}_{l_1}:\pi_1(U\setminus (U\cap \overline B),p_1)\to
\Sigma_d$, and the paths $l_{o_i}$ define homomorphisms
$\text{im}_{l_{o_i}}:\pi_1(U_{o_i}\setminus (U_{o_i}\cap \overline
B),p_{oi})\to \pi_1(U\setminus (U\cap \overline B),p_1)$ and
$\psi_{l_{o_i}}=f_*\circ \text{im}_{l_{o_i}}:\pi_1(U_{o_i}\setminus
(U_{o_i}\cap \overline B),p_{o_i})\to \Sigma_d$. Denote by
$H_{\overline B_1}$ the image $\varphi(\pi_1(U\setminus
(U\cap\overline B),p_{1}))$ in $G$. It is easy to see that
$\varphi_{\widetilde 1_{o_i}}=\psi_{l_{o_i}}$. Therefore,
$H_{B_1}\subset H_{\overline B_1}$ and $H_{o_I}\subset H_{\overline
B_1}$.

Let $\gamma_{q_{o_i}}\in \pi_1(U_{o_i}\setminus (U_{o_i}\cap
\overline B),p_{o_i})$ be a geometric generator represented by the
circle $\delta_{o_i}=\{ q_{o_i}\}\times \{ u_1\in \mathbb C \mid
|u_1|= \frac{1}{2}\}\subset \{ q_{o_i}\}\times D_1\subset
U_{o_i}\setminus \overline B$ (the circuit in positive direction).
It is easy to see that $\text{im}_{l_{o_i}}(\gamma_{q_{o_i}})=
\gamma_{q_{1}}$. Therefore $\varphi_{\widetilde
1_{o_i}}(\gamma_{q_{o_i}})=\sigma_1$.

Consider the restriction of $f$ to each $V_{o_i,m}$,
$f_{i,m}=f_{|V_{o_i,m}}:V_{o_i,m}\to U_{o_i}$. Denote by $d_{i,m}$
the degree of $f_{i,m}$, $d=d_{i,1}+\dots +d_{i,k_i}$. By
construction, for the point $\overline o_{i,m}=V_{o_i,m}\cap
f^{-1}(o_i)$ we have $\deg_{\overline o_{i,m}}f=d_{i,m}$.

The numbering of the points of $f^{-1}(p_0)$ and the path
$\widetilde l_{o_i}$ define a numbering of the points of
$f^{-1}(p_{o_i})$. Then the decomposition
$f^{-1}(U_{o_i})=V_{o_i,1}\sqcup\dots \sqcup V_{o_i,k_i}$ defines a
partition $J_i=\{ I_{i,1},\dots , I_{i,k_i}\}$ of $[1,d]$, $j\in
I_{i,m}$ if and only if $\widetilde p_j\in f^{-1}(p_{o_i})\cap
V_{o_i,m}$. By construction, the group $H_{o_i}$ leaves invariant
the partition $J_i$ and acts transitively on each $I_{i,m}\in J_i$.
In particular, the action of $\sigma_1$ leaves invariant the
partition $J_i$.

Assume that if $\overline B_j$ is an irreducible component of the
branch locus $\overline B$ of a covering $f$ such that
$f_*(\gamma_{q_j})$ is not a transposition, then $f$ satisfies the
following conditions:
\begin{itemize}
\item[($R_1)$]
{\it For each $o_i\in \overline B_j\cap\text{Sing}\, \overline B$
and each $V_{o_i,m}$ there is at most one irreducible component
$\overline R_k\subset f^{-1}(\overline B_j)$ of the ramification
locus of $f$ which intersects with $V_{o_i,m}$.
\item[($R_2$)] For each $\overline R_k\subset f^{-1}(\overline B_j)$
there is $o_i\in \overline B_j\cap\text{Sing}\, \overline B$ and $m$
such that $\overline R_k\cap V_{o_i,m}\neq \emptyset$ and
$r_k<d_{i,m}$.
\item[($R_3$)] If $\overline R_k\cap V_{o_i,m}\neq \emptyset$
and $\widetilde R$ is another ramification curve of $f$ such that
$\widetilde R\cap V_{o_i,m}\neq \emptyset$, then for a point $q\in
f(\widetilde R)$ the image $f_*(\gamma_q)$ of a geometric generator
$\gamma_q$ is a transposition in} $\Sigma_d$.
\end{itemize}

\begin{lem}\label{H+J1} Let $f$ and its branch curve $\overline B_1$
satisfy conditions ($R_0$) -- ($R_3$), and let $H$ be a subgroup of
$H_{\overline B_1}$ generated by $\sigma_1$ and the transpositions
belonging to $H_{\overline B_1}$. Then $H=\Sigma_{J_{H}}$.
\end{lem}
\proof Let $\sigma=\sigma_1=\varphi_1( \gamma_{q_1})$ where
$\gamma_{p_1}$ is the geometric generator defined above. Then it is
easy to see that 
condition ($R_1$) implies that $H$ and $\sigma$ satisfy condition
($i$) from Lemma \ref{H+J}. Similarly, it follows from conditions
($R_2$) and ($R_3$) that $H$ and $\sigma$ satisfy conditions ($ii$)
and ($iii$) from Lemma \ref{H+J}. Therefore, $H=\Sigma_{J_{H}}$.
\qed

\begin{prop}\label{G+J1} Let $Z_1,\dots,Z_k$ be the irreducible
components $f:Z\to Y$ be a ramified covering of a simply connected
surface $Y$. Assume that the branch locus $\overline B$ of $f$
satisfies conditions ($R_0$) -- ($R_3$). Then the Galois group $G$
of $f$ is isomorphic to $\Sigma_{d_1}\times \dots \times
\Sigma_{d_k}$, where  $d_i=\deg f_{|Z_i}$.
\end{prop}
\proof The decomposition $Z=Z_1\sqcup\dots\sqcup Z_k$ defines a
partition $J=\{ I_1,\dots, I_k\}$ of the set $f^{-1}(p_0)$. The
group $G$ leaves invariant the partition $J$ and acts transitively
on each $I_j\subset J$. Therefore Proposition \ref{G+J1} follows
from Lemma \ref{H+J1} and Corollary \ref{coro1}. \qed

\section{Proof of Theorem \ref{main0}, \ref{main1} and Corollary \ref{main2}}\label{thm1}
We use notations defined in Introduction and Section
\ref{cover}.
\subsection{Proof of Theorem \ref{main1}}
Denote by $\overline C_i=\nu^{-1}(C_i)$ the irreducible components
of $\overline C$, $1\leq i\leq k$.

Obviously,
$$Z\simeq \{  (p,l)\in\overline C\times \hat{\mathbb P}^2\mid
p\in \overline C,\, \nu(p)\in l\} $$
and it is easy to see that
$$Z_i\simeq \{  (p,l)\in\overline C_i\times \hat{\mathbb P}^2\mid
p\in \overline C_i,\, \nu(p)\in l\} $$
are the irreducible components of the surface $Z$.

Let $t$ be a local parameter in a small neighbourhood $U\subset
\overline C$ of a point $p\in \overline C$ and let the normalization
$\nu$ be given in $U$ by
\begin{equation} \label{nu1} x_1=\phi_1(t), \, \,
x_2=\phi_2(t), \, \, x_3=\phi_3(t). \end{equation} If
$(y_1,y_2,y_3)$ are homogeneous coordinates in $\hat{\mathbb P}^2$
dual to the coordinates $(x_1,x_2,x_3)$ in $\mathbb P^2$, then the
surface $Z$, in the neighbourhood $U\times \hat{\mathbb P}^2\subset
\overline C\times \hat{\mathbb P}^2$, is given by equation
$$y_1\phi_1(t)+y_2\phi_2(t)+y_3\phi_3(t)=0.$$
In particular, if $\nu$ is given by equations (\ref{nu}), that is,
\begin{equation} \label{nu5} \phi_1=\sum_{i=s_p}^{\infty} a_it^i, \, \,
\phi_2=t^{r_p}, \, \, \phi_3=1,
\end{equation}
then $Z\cap (U\times \hat{\mathbb P}^2)$ lies in $U\times {\mathbb
C}^2$, where $\mathbb C^2=\{ y_1\neq 0\}$ is the affine plane in
$\hat{\mathbb P}^2$, and $Z\cap (U\times {\mathbb C}^2)$ is given by
equation
\begin{equation} \label{nu2}
\sum_{i=s_p}^{\infty} a_it^i +z_2t^{r_p}+z_3=0,
\end{equation} where $z_i=y_i/y_1$, $i=2,3$. Therefore $Z$ is a
smooth surface and $(t,z_2)$ are coordinates in $Z\cap (U\times
{\mathbb C}^2)$.

The restriction of the covering $f$ to $Z\cap (U\times {\mathbb C}^2)$,
$$f_U:Z\cap (U\times {\mathbb C}^2)\to \mathbb C^2,$$
is the restriction of the projection $(t,z_2,z_3)\mapsto (z_2,z_3)$,
therefore it is given by
\begin{equation}\label{covf}
\begin{array}{ll}
z_2 &= z_2, \\
z_3 &= -(\displaystyle \sum_{i=s_p}^{\infty} a_it^i +z_2t^{r_p}).
\end{array}\end{equation}
Its Jacobian is equal
$$J(f_U)=-t^{r_p-1}(\sum_{i=s_p}^{\infty}i a_it^{i-r_p} +r_pz_2).$$
Therefore $f_U$ is ramified along a curve $R$ given by equation
\begin{equation}\label{ram}
\frac{1}{r_p}\sum_{i=s_p}^{\infty} a_it^{i-r_p} +z_2=0
\end{equation}
with multiplicity two and along the fibre $F_p=\{ t=0\}$ with
multiplicity $r_p$ if $r_p\geq 2$ and hence $f(F_p)=L_{p}\subset
\hat{\mathbb P}^2$ is a component of the branch locus of $f$ if
$r_p\geq 2$. Note also that $R$ is a section of the ruled surface
$Z\cap (U\times {\mathbb C}^2)\to \mathbb C^2$ and, in addition, it
is the unique section contained in the ramification locus. Therefore
to show that $R$ is a germ of the curve $\widetilde C$, we can
assume that the image $\nu(p)$ is a smooth point of $C$, that is, we
can assume that $r_p=1$ and $\phi_2(t)= t$. Then $R$ is given by
$\phi'_1(t)+z_2=0$ and the restriction of $f_U$ to $R$ is given by
\begin{equation} \label{dual} y_1=1,\, \, y_2=-\phi_1'(t),\, \,
y_3=-\phi_1(t)+t\phi_1'(t).\end{equation} Everyone easily check that
that equations (\ref{dual}) together with the equations
$$ x_1=\phi_1(t),\, \, x_2=t,\, \, x_3=1$$
(defining the germ $\nu(U)$ of $C$) is a parametrization of $\check
C\subset I$ over $\nu(U)$.

To count the local degree $\deg_{q}f_U$ of the covering $f_U:Z\cap
(U\times {\mathbb C}^2)\to \mathbb C^2$ at the point $q=(0,0,0)$,
first of all, note that the curve $\{ z_2=0\}$ does not belong to
the ramification locus of $f_U$, since $a_{s_p}\neq 0$ in equation
(\ref{ram}). Next, let us choose a new parameter $t_1$ such that
$t_1^{s_p}=
\sum_{i=s_p}^{\infty} a_it^{i}$, then $f_U$
is given by equations of the form
\begin{equation}\label{covf1}
\begin{array}{ll}
z_2 &= z_2, \\
z_3 &= -(t_1^{s_p}+ z_2
\sum_{i=r_p}^{\infty} c_it_1^i
)
\end{array}\end{equation}
and to count $\deg_{q}f_U$, it suffices to count the number of
points belonging to $f_U^{-1}((z_{2,0},z_{3,0}))$, where a point
$(z_{2,0},z_{3,0})\in Im f_U$ is such that $z_{2,0}=0$, $z_{3,0}\neq
0$, and $z_{3,0}$ is sufficiently close to zero. It follows from
equations (\ref{covf1}) that this number is equal to $s_p$. \qed

\subsection{Proof of Theorem \ref{main0}}
By Theorem \ref{main1}, the branch locus $\overline B$ of the
dualizing covering $f:Z\to\hat{\mathbb P}^2$ consists of the curve
$\hat C$ and the lines $L_p$ for which $r_p\geq 2$, the ramification
locus $\overline R$ consists of the curve $\widetilde  C$ and the
fibres $F_p$ with $r_p\geq 2$. Each line $L_p$ and each fibre $F_p$
are smooth and $F_{p_1}\cap F_{p_2}=\emptyset$ for $p_1\neq p_2$.
Next, the restriction $f_{|F_p}: F_p \to L_p$ of $f$ to $F_p$ is an
isomorphism and the restriction $f_{|\widetilde C}: \widetilde C\to
\hat C$ of $f$ to $\widetilde C$ is a bi-rational map, and $f$ is
ramified along $\widetilde C$ with multiplicity two. Therefore
$f_*(\gamma_l)$ is a transposition for any geometric generator
$\gamma_l\in \pi_1(\hat{\mathbb P}^2\setminus \overline B,l_0)$,
$l\in \hat C$, and hence the dualizing covering $f$ and its branch
locus $\overline B$ satisfy conditions ($R_0$) and ($R_1$) (see
Section \ref{cover}). Next, if $F_p\subset \overline R$ then the
point $q=F_p\cap \widetilde  C$ belongs to $f^{-1}(\text{Sing}\,
\overline B)$ and $s_p=\deg_q f>r_p$, that is, $f$ and its branch
locus $\overline B$ satisfy conditions ($R_2$) and ($R_3$). Now
Theorem \ref{main0} follows from Proposition \ref{G+J1}. \qed

\subsection{Proof of Corollary \ref{main2}}
First of all, note that if a line $l\subset \mathbb P^2$ is an
irreducible component of the curve $C$ then for each loop
$l_t\subset T_C$ starting at $l_0$ the monodromy defined by $l_t$
leaves fixed the point $l\cap l_0$.

Let $C_1,\dots ,C_n$ be the irreducible components of $C$ and let
$\deg C_i=d_i\geq 2$ for $i=1,\dots, k$ and  $\deg C_i=1$ for $i>k$.
Denote by $\text{Sing}\, C$ the set of singular points of $C$, by
$\displaystyle L_{Sing}=\bigcup_{P\in \text{Sing}\, C}L_P\subset
\hat{\mathbb P}^2$, by $C'=C_1\cup\dots\cup C_k$. Let $f':Z\to
\hat{\mathbb P}^2$ be the dualizing covering for $C'$ and $\overline
B'$ its branch locus (see Theorem \ref{main1}). Then it is easy to
see that $T_C=\hat{\mathbb P}^2\setminus (\overline{B}'\cup
L_{Sing})$.

We have $f'^{-1}(l_0)=\{ (p_1,l_{0}),(p_2,l_{0}),\dots, (p_d,l_0)
\}$, where $p_i=\nu^{-1}(P_i)$, $l_0\cap C'=\{P_1,P_2,\dots ,P_d\}$,
and $P_1,P_2\in C_1$.

The embedding $i:T_C\hookrightarrow \hat{\mathbb P}^2\setminus
\overline B'$ defines an epimorphism $i_*: \pi_1(T_C,l_0)\to
\pi_1(\hat{\mathbb P}^2\setminus \overline B',l_0)$. By Theorem
\ref{main0}, there is an element $\gamma'\in \pi_1(\hat{\mathbb
P}^2\setminus \overline B',l_0)$ such that $f'_*(\gamma') =(1,2)\in
\Sigma_d$, where $(1,2)$ is the transposition transposing the points
$(p_1,l_{0})$ and $(p_2,l_{0})$. Let $l_t$ be a loop representing an
element $\gamma \in \pi_1(T_C,l_0)$ such that $i_*(\gamma)=\gamma'$.
It is easy to see that the loop $l_t$  is a desired one. \qed

 \ifx\undefined\bysame
\newcommand{\bysame}{\leavevmode\hbox to3em{\hrulefill}\,}
\fi

\end{document}